\documentclass[12pt, a4paper, reqno]{amsart}

\usepackage{enumerate}
\usepackage{amssymb}

\newtheorem{lemma}{Lemma}[section]

\newtheorem{theorem}[lemma]{Theorem}
\newtheorem{proposition}[lemma]{Proposition}
\newtheorem{remark}[lemma]{Remark}

 \newcommand{\C}{\mathbb{C}}
 \DeclareMathOperator{\Div}{div}

\author[A. Gasull]{A. Gasull$^*$}
\address{Dept. de Matem\`{a}tiques.
Universitat Aut\`{o}noma de Barcelona. Edifici C. 08193 Bellaterra, Barcelona.
Spain} \email{gasull@mat.uab.es}

\author{H. Giacomini}
\address{Laboratoire de Math\'{e}matiques et Physique Th\'{e}orique. Facult\'{e} des
Sciences et Techniques. Universit\'{e} de Tours. 37200 Tours. France}
\email{Hector.Giacomini@phys.univ-tours.fr}

\author[J. Torregrosa]{J. Torregrosa$^*$}
\address{Dept. de Matem\`{a}tiques \\
Universitat Aut\`{o}noma de Barcelona \\ Edifici C. 08193 Bellaterra, Barcelona.
Spain} \email{torre@mat.uab.es}

\thanks{$^*$\uppercase{S}upported by \uppercase{DGES
N}o.\uppercase{BFM}2002-04236-\uppercase{C}02-2
 and \uppercase{CONACIT} 2001\uppercase{SGR}-00173.}

\date{}
\dedicatory{} \commby{}

\begin{document}

\title[Explicit non-algebraic limit cycles]{Explicit non-algebraic limit cycles
for polynomial systems}

\begin{abstract}
We consider a system of the form $\dot x\, =\, P_n(x,y)+xR_m(x,y)$, $\dot y\,=\, Q_n(x,y)+yR_m(x,y)$, where $P_n(x,y)$, $Q_n(x,y)$ and $R_m(x,y)$ are homogeneous
polynomials of degrees $n$, $n$ and $m$, respectively, with $n\le
m$. We prove that this system has at most one limit cycle and that when it exists it can be explicitly found. Then we study a particular case, with $n=3$ and $m=4$. We prove that this quintic polynomial system has an explicit limit cycle which is not algebraic. To our knowledge,
there are no such type of examples in the literature.\\
The method that we introduce to prove that this limit cycle is not algebraic
can be also used to detect algebraic solutions for other families
of polynomial vector fields or for probing the absence of such type of solutions.
\end{abstract}

\maketitle

\section{Introduction and Main Results}

Examples of planar polynomial vector fields having explicit
algebraic limit cycles appear in most textbooks of ordinary
differential equations. One of the simplest examples is the one of
a cubic system that in polar coordinates writes as $\dot
r=r(1-r^2), $ $\dot \theta=1,$ see for instance \cite{Per2001}. On
the other hand, it seems intuitively clear that ``most'' limit
cycles of planar polynomial vector fields have to be
non-algebraic. Nevertheless, until 1995 it was not proved that the
limit cycle of the van der Pol equation is not algebraic, see
\cite{Oda1995}.

The goal of this paper is to give a planar polynomial vector field
for which we can get an explicit limit cycle which is not
algebraic. As far as we know there are no examples of this
situation in the literature.

Recall that a real or complex polynomial
$F(x,y)$ is an \textit{algebraic solution} of a real
polynomial system $(\dot x, \dot y)=(X(x,y), Y(x,y))$
if
\begin{equation}\label{eq:1}
\frac{\partial F(x,y)}{\partial x}X(x,y)+ \frac{\partial
F(x,y)}{\partial y}Y(x,y)= K(x,y)F(x,y),
\end{equation}
for some polynomial $K(x,y),$ called the \textit{cofactor} of $F.$
Notice that when $F(x,y)$ is real, the curve $F(x,y)=0$ is
invariant under the flow of the differential equation. Observe also
that the degree of the cofactor is one less than the degree of the
vector field. A limit cycle is called algebraic if it is an oval
of a real algebraic solution.

Our main result is:

\begin{theorem}\label{th:1.1}
The planar differential system
\begin{equation}\label{eq:2}
\left\{
\begin{array}{ccl}
\dot x&=&-(x-y)(x^2-xy+y^2)+x(2x^4+2x^2y^2+y^4), \\
\dot y&=&-(x+y)(2x^2-xy+2y^2)+y(2x^4+2x^2y^2+y^4),
\end{array}
\right.
\end{equation}
has exactly one limit cycle which is hyperbolic and non-algebraic.
In polar coordinates, this limit cycle is
\[
r=e^{\frac{3}{2}\psi(\theta)-\theta}
\left(a+ 2\int_0^\theta{\frac{\cos^4(s)+1}{\cos^2(s)+1}
e^{3\psi(s)- 2s}}ds\right)^{-1/2},
\]
where $\psi(\theta)=2\int_0^\theta \frac{1+\tan^2s}{2+\tan^2s}ds$
and
\[
a=\frac{2e^{4\pi-6\sqrt{2}\pi}}{1-e^{4\pi-6\sqrt{2}\pi}}
\int_0^{2\pi}\frac{\cos^4(s)+1}{\cos^2(s)+1}
e^{3\psi(s)-2s}ds\approx 1.19903.
\]
\end{theorem}

The main steps of the proof are:
\begin{enumerate}[(i)]
\item We consider the family of systems,
\begin{equation}\label{eq:3}
\left\{
\begin{array}{ccl}
\dot x&=&P_n(x,y)+xR_m(x,y), \\
\dot y&=&Q_n(x,y)+yR_m(x,y),
\end{array}
\right.
\end{equation}
where $P_n(x,y)$, $Q_n(x,y)$ and $R_m(x,y)$ are homogeneous
polynomials of degrees $n,$ $n$ and $m,$ respectively, with $n\le
m.$ In Section~\ref{se:2} we prove that all of them have at most one
limit cycle and that when it exists it can be explicitly found.

\item We then consider a concrete system of the form~(\ref{eq:3}) with $n=3$ and $m=4,$
having an explicit, unique and hyperbolic limit cycle, see
Section~\ref{se:4}.

\item Finally, we use the method developed in Section~\ref{se:5} for
studying the algebraic solutions of the fixed system. Once all of them are found, none of which has ovals,
we can prove that the limit cycle is non-algebraic, as shown in Section~\ref{se:6}.
\end{enumerate}

In Section~\ref{se:3} we study system~(\ref{eq:3}) with $n=1.$ In
this case we prove that when the limit cycle exists it is always
algebraic. This is a short preliminary that we include for the
sake of completeness and that helps to understand how we have
arrived to the final example studied in our main theorem.

We end this introduction by showing a system of the form~(\ref{eq:3}), almost equal to the one studied in Theorem~\ref{th:1.1}, having also a unique hyperbolic limit cycle, but
algebraic. Concretely, it is easy to check that the system
\begin{equation}\label{eq:4}
\left\{
\begin{array}{ccl}
\dot x&=&-(x-y)(x^2-xy+y^2)+x(x^4+3x^2y^2+2y^4), \\
\dot y&=&-(x+y)(2x^2-xy+2y^2)+y(x^4+3x^2y^2+2y^4)
\end{array}
\right.
\end{equation}
has the algebraic limit cycle $1-x^2-y^2=0.$
The results when $n=1$ and systems~(\ref{eq:2}) and (\ref{eq:4}) illustrate that the problem of distinguishing whether the limit cycle for system~(\ref{eq:3}) is algebraic or
not is quite hidden in the coefficients of the system.

\section{Systems with explicit limit cycles}\label{se:2}

This section is devoted to study system~(\ref{eq:3}). In polar
coordinates it writes as
\begin{equation}\label{eq:5}
\left\{
\begin{array}{ccl}
\dot r&=&f(\theta)r^n+h(\theta)r^{m+1}, \\
\dot \theta&=&g(\theta)r^{n-1},
\end{array}
\right.
\end{equation}
where
\[
\begin{aligned}
f(\theta)&=\cos \theta P_n(\cos \theta, \sin \theta)+
\sin \theta Q_n(\cos \theta , \sin \theta ) , \\
g(\theta)&=\cos \theta Q_n(\cos \theta , \sin \theta )-
\sin \theta P_n(\cos \theta , \sin \theta ), \\
h(\theta)&=R_m(\cos \theta, \sin\theta ).\\
\end{aligned}
\]

\begin{theorem}\label{th:2.1}
System~(\ref{eq:3}) has at most one limit cycle. When it exists
it is hyperbolic and in polar coordinates it writes as
\[
r=\left(\exp\left[\int_{0}^{\theta}\frac{f(s)}{g(s)}ds\right]
\left [a + \int_{0}^{\theta}\frac{h(s)}{g(s)}\exp\left(-\int
_{0}^s\frac{f(w)}{g(w)}dw\right)ds\right]\right)^{\frac{1}{n-m-1}},
\]
where $a=AB/(1-A),$ being
\[
A=\exp\left(\int_{0}^{2\pi}\frac{f(s)}{g(s)}ds\right) \quad\mbox{ and
}\quad B=\int_{0}^ {2\pi}\frac{h(s)}{g(s)}\exp\left(-\int
_{0}^s\frac{f(w)}{g(w)}dw\right)ds.
\]
\end{theorem}

\begin{proof}
Consider the expression of system~(\ref{eq:3}) in polar
coordinates, \textit{i.e. } system~(\ref{eq:5}). If $g(\theta)$
vanishes for some $\theta=\theta^*$ then it has
$\{\theta=\theta^*\}$ as an invariant straight line. From the uniqueness of solutions we get that system~(\ref{eq:3}) has no limit cycles. If $g(\theta)\ne 0$ then we can write the
system as
\[
\frac{dr}{d\theta}=\frac{f(\theta)}{g(\theta)}r+\frac{h(\theta)}{g(\theta)}r^{m-
n + 2},
\]
which is a Bernoulli equation. By introducing the standard change of variables
$\rho=r^{n-m-1}$ we obtain the linear equation
\begin{equation}\label{eq:6}
\frac{d\rho}{d\theta}=(n-m-1)\frac{f(\theta)}{g(\theta)}\rho+(n-m-1)\frac{h(\theta)}{g(\theta)}.
\end{equation}
Notice that system~(\ref{eq:3}) has a periodic orbit if and only if
equation~(\ref{eq:6}) has a strictly positive $2\pi$ periodic solution.

The general solution of equation~(\ref{eq:6}), with
initial condition $\rho(0)=\rho_0$, is
\begin{equation}\label{eq:7}
\rho(\theta;\rho_0)=\exp\left(\int_{0}^{\theta}\frac{f(s)}{g(s)}ds\right)\!\!\left(\rho_0+\int_{0}
^{\theta}\frac{h(s)}{g(s)}\exp\left(-\int_{0}^s\frac{f(w)}{g(w)}dw\right)ds\right)\!:=G(\theta, \rho_0).
\end{equation}
The condition that the solution starting
at $\rho=\rho_0$ is periodic reads as $A(\rho_0+B)=\rho_0$. Hence, if $A=1$ and
$B=0$, system~(\ref{eq:3}) has a continuum of periodic orbits,
otherwise it has at most the solution starting at $\rho_0=AB/(1-A):=a$.

In order to prove the hyperbolicity of the limit cycle
notice that the Poincar\'e return map is $\Pi(\rho_0)=\rho(2\pi;\rho_0)$.
Thus $\Pi'(\rho_0)=\exp\left(\int_{0}^{2\pi} \frac{f(s)}{g(s)}ds\right)=A\ne 1$
for all $\rho_0$, and in particular we get that the limit cycle is hyperbolic, whenever it exists.
\end{proof}

>From the proof of the above theorem we also get the following
remark.

\begin{remark}\label{re:2.2} When system~(\ref{eq:3}) has a limit cycle
it can be written in the form $F(r, \theta):=r^{m-n+1}G(\theta,
a)-1=0.$ As we will see, the function $F(r, \theta)$ can be
algebraic or not in cartesian coordinates, depending on the
concrete system considered. In any case the expression given in
(\ref{eq:1}) can be also extended to non-algebraic functions $F$,
and in this case the cofactor $K$ is not necessarily a polynomial.
Curiously enough, independently of the algebraiticity of $F$ its
corresponding cofactor satisfying~(\ref{eq:1}) is always the
polynomial $K(x,y)=(m-n+1)R_m(x,y).$ Examples of non-algebraic
solutions having a polynomial cofactor have been already given in
the literature, see for instance \cite{GarGin2003},
\cite{GiaGinGra2005}. The example presented in
Theorem~\ref{th:1.1} provides a non-algebraic limit cycle having a
polynomial cofactor.
\end{remark}

The next remark explains why our proof that the limit cycle of system~(\ref{eq:2}) is not algebraic
does not use the explicit expression
of the limit cycle. On the contrary, in Section~\ref{se:5}, we
develop a method for studying all the algebraic solutions of a
system having at least a solution of the form $y-\alpha(x)=0,$ where $\alpha(x)$ is a rational
function, and we apply this method to determine all the algebraic solutions of system~(\ref{eq:2}).

\begin{remark}\label{re:2.3}
Although for system~(\ref{eq:3}) we know that the expression of
the limit cycle is $r^{m-n+1}G(\theta, a)-1=0$, it is not an easy
task to elucidate whether this curve is algebraic or not in cartesian coordinates. As a
example of this difficulty we recall the Filiptsov's example, see
\cite{Fil1973}:
\[
\left\{
\begin{array}{ccl}
\dot x&=&6(1+a)x+2y-6(2+a)x^2+12xy, \\
\dot y&=&15(1+a)y+3a(1+a)x^2-2(9+5a)xy+16y^2,
\end{array}
\right.
\]
which has the algebraic solution
$3(1+a)(ax^2+y^2)+2y^2(2y-3(1+a)x)$. This algebraic solution contains a limit cycle for $0<a<3/13$. For the sake of simplicity we fix $a=1/6.$ For this value the limit cycle
is $rG(\theta)-1=0$, where
\[
G(\theta)=\frac{7(\sin^4\theta-2\sin^2\theta+1)}
{6\sin\theta(-17\sin^2\theta+42\sin\theta\cos\theta-7\pm
2\sqrt{\varphi(\theta)})}
\]
and
\[
\varphi(\theta)=60\sin^4\theta-357\sin^3\theta\cos\theta
+84\sin^2\theta+441\sin^2\theta\cos^2\theta-147\sin\theta\cos\theta.
\]
Note that is not easy at all to realize that the expression
$rG(\theta)-1=0$ corresponds to a polynomial in cartesian coordinates.
\end{remark}

It is not difficult to see that system~(\ref{eq:3}) has always
algebraic solutions.
\begin{lemma}\label{le:2.4}
System~(\ref{eq:3}) has $F(x,y)=yP_n(x,y)-xQ_n(x,y)$ as an
algebraic solution with cofactor $(n+1)R_m+\Div(P_n, Q_n).$ Notice
that it is formed by a product of (complex or real) invariant
straight lines through the origin.
\end{lemma}
\begin{proof}
By using the homogeneity of $P_n$ and $Q_n$, we know from the
Euler's formula that $nP_n(x,y)=x\frac{\partial P_n}{\partial
x}+y\frac{\partial P_n}{\partial y}$ and
$nQ_n(x,y)=x\frac{\partial Q_n}{\partial x}+y\frac{\partial
Q_n}{\partial y}.$ Thus
\begin{eqnarray*}
\left(y\frac{\partial P_n}{\partial x}-Q_n-x\frac{\partial Q_n}{\partial x}\right)
(P_n+xR_m)+\left(P_n+y\frac{\partial P_n}{\partial y}-x\frac{\partial Q_n}{\partial y}\right)(Q_n+yR_m)\\
=\left((n+1)R_m+\frac{\partial P_n}{\partial x}+\frac{\partial Q_n}{\partial y}\right)F.
\end{eqnarray*}
\end{proof}

\begin{lemma}\label{le:2.5}
Let $F(x,y)$ be an algebraic solution of degree $\ell$ of the
system
\[
\left\{
\begin{array}{ccl}
\dot x&=&P(x,y)+xR_m(x,y), \\
\dot y&=&Q(x,y)+yR_m(x,y),
\end{array}
\right.
\]
where $P(x,y)$ and $Q(x,y)$ are polynomials of degree less or
equal than $n$ and $R_m(x,y)$ is a homogeneous polynomial of
degree $m$, with $n\le m.$ Thus the homogeneous part of maximum
degree of its cofactor is $\ell R_m(x,y)$.
\end{lemma}
\begin{proof}
Since $F$ is an algebraic solution of system~(\ref{eq:3}) we know that
\[
\frac{\partial F}{\partial x}\left(P+xR_m\right) +
\frac{\partial F}{\partial y}\left(Q+yR_m\right)= KF,
\]
where $K$ is the cofactor of $F.$

Denote by $F_\ell(x,y)$ and by $K_m(x,y)$ the homogeneous parts
of maximum degree of $F(x,y)$ and $K(x,y),$ respectively. By using
the homogeneity of $F_\ell$ we know, from the Euler's formula,
that $\ell F_\ell=x\frac{\partial F_\ell}{\partial x}+y\frac{\partial F_\ell}{\partial y}$. By
equating the higher degree terms in the above equation we obtain
\[
 \ell F_\ell R_m=\frac{\partial F_\ell}{\partial x}xR_m+\frac{\partial
F_\ell}{\partial y}yR_m = K_m F_\ell.
\]
Thus $K_m(x,y)=\ell R_m(x,y)$ as we wanted to prove.
\end{proof}

Next lemma collects some easy remarks on the structure of the
cofactors.

\begin{lemma}\label{le:2.6}
Let
\[
\left\{
\begin{array}{ccl}
\dot x&=&X(x,y), \\
\dot y&=&Y(x,y),
\end{array}
\right.
\]
be a real planar polynomial system. The following holds:
\begin{enumerate}[(i)]
\item If it has a complex algebraic solution, then it also has a
real algebraic solution.

\item Assume the vector field satisfies
\[
(X(-x, -y), Y(-x, -y))=(-1)^s(X(x,y), Y(x,y)), \quad \mbox{being $s$
either 0 or 1.}
\]
 Then if it has a
real algebraic solution then it has another real algebraic
solution with cofactor $K$ satisfying $K(-x, -y)=(-1)^{s+1}K(x,y).$
\end{enumerate}
\end{lemma}
\begin{proof} The proof is elementary.

%
%
\end{proof}

Finally we give an integrating factor for system~(\ref{eq:3}).

\begin{lemma}\label{le:2.7} Consider system~(\ref{eq:3}) and define
\[
V(x,y)=(r^{m-n+1}G(\theta, a)-1)(yP_n(x,y)-xQ_n(x,y)),
\]
where $G(\theta, \rho_0)$ is
the function given in (\ref{eq:7})
and $\rho_0=a$ is the value for which this function is
$2\pi-$periodic. Then, whenever it is defined, $1/V(x,y)$ is an
integrating factor of the system and we call $V(x,y)$ an inverse integrating factor.
\end{lemma}

\begin{proof}
We use the following formula: let $F_1$ and $F_2$ be two solutions of
$(\dot x, \dot y)=(X(x,y), Y(x,y))$ with cofactors $K_1$ and $K_2,$
respectively. Thus
\[
\Div\left(\frac{(X, Y)}{F_1F_2}\right)=\frac1{F_1F_2}\left(\Div(X, Y)-(K_1+K_2)\right).
\]
We remark that the above formula, taking a denominator of the form
$\prod F_i^{\alpha_i},$ for some real or complex constants
$\alpha_i,$ is indeed the key point of the Darboux theory of
integrability, see \cite{Lli2004}. Take
$F_1(x,y)=yP_n(x,y)-xQ_n(x,y)$ and $F_2(x,y)=r^{m-n+1}G(\theta,
a)-1.$ By using Lemma~\ref{le:2.4} and Remark~\ref{re:2.2} we know
that their associated cofactors areit being non-Liouvillian
$K_1(x,y)=(n+1)R_m(x,y)+\Div(P_n(x,y), Q_n(x,y))$ and
$K_2(x,y)=(m-n+1)R_m(x,y),$ respectively. On the other hand,
taking the vector field associated to system~(\ref{eq:3}) we get
\[
\Div(X, Y)=\Div(P_n, Q_n)+2R_m+x\frac{\partial R_m}{\partial x}+
y\frac{\partial R_m}{\partial y}=\Div(P_n, Q_n)+(2+m)R_m,
\]
where we have used again Euler's formula. Collecting all the above results
we get $\Div\left((X, Y)/(F_1F_2)\right)\equiv0$ as
we wanted to prove.
\end{proof}

\begin{remark} (i) When we apply the above Lemma to systems~(\ref{eq:2}) and (\ref{eq:4}) we get both non-algebraic and
algebraic inverse integrating factors.

(ii) In \cite{GiaLliVia1996} it is proved that when $1/V(x,y)$ is
an integrating factor of $(\dot x, \dot y)=(X(x,y), Y(x,y))$ and
$V(x,y)$ is defined in the whole plane, all the limit cycles of
the system are included in the curve $V(x,y)=0.$ This is the case
of system~(\ref{eq:3}); the limit cycle, whenever it exists, it is
given by the expression $F_2(x,y)=r^{m-n+1}G(\theta, a)-1=0.$

(iii) The equality $\Div\left((X, Y)/(F_1F_2)\right)\equiv0$
also holds when instead of
$F_2(x,y)=r^{m-n+1}G(\theta, \rho_0)|_{\rho_0=a}-1$ we take a
different value of $\rho_0,$ but in this case $F_2$ is indeed a
multivaluated function and the result of \cite{GiaLliVia1996} can not be applied.
\end{remark}

\section{A family of systems with explicit algebraic limit cycles}\label{se:3}
The existence of limit cycles for a subfamily of system~(\ref{eq:3}) has been studied
in \cite{GasTor2004}. Here we prove that the limit cycle found there is algebraic.

\begin{proposition}\label{pr:3.1}
Consider the system
\begin{equation}\label{eq:8}
\left\{
\begin{array}{ccl}
\dot x&=&-y+x(a+R_m(x,y)), \\
\dot y&=&\phantom{-}x+y(a+R_m(x,y)),
\end{array}
\right.
\end{equation}
where $a$ is a real parameter and $R_m(x,y)$ is a homogeneous
polynomial of degree $m.$ Then it has only two algebraic invariant
curves $x^2+y^2$ and $H(x,y)=G_m(x,y)-1,$ where
$G(\theta)=G_m(\cos\theta, \sin\theta)$ satisfies
$G'+maG+mR_m(\cos\theta, \sin\theta)=0$. Furthermore, when the limit
cycle exists, $m$ is even and $H(x,y)$ contains a real oval which
is the limit cycle of the system.
\end{proposition}

\begin{proof}
>From Lemma~\ref{le:2.4} we can see that $x^2+y^2$ is an
algebraic solution with cofactor $K(x,y)=2a+2R_m(x,y)$. Now, we study other possible algebraic
solutions.

Write the Fourier expansion of $R_m$:
\[
R_m(\cos\theta, \sin\theta)=\sum_{k=-m}^{k=m} c_{k}e^{k\theta i},
\mbox{ where }\overline{c_k}=c_ {-k}\in \C, c_k=0 \mbox{ when } k \not\equiv m (\textrm{mod }2).
\]
Following the steps of the proof of Theorem~\ref{th:2.1} we
obtain that, in polar coordinates, the solution of (\ref{eq:8})
starting at $r=r_0$ when $\theta=0$ can be written as
\[
r^{-m} = \left(r_0^{-m}+m\sum\limits_{k=-m}^{k=m} \frac{c_{k}}
{ki+ma}\right)e^{-ma\theta}+G_m(\cos \theta , \sin \theta ),
\]
or
\[
1 = \left(r_0^{-m}+m\sum\limits_{k=-m}^{k=m}\frac{c_{k}}
{ki+ma}\right)r^me^{-ma\theta} + G_m(r\cos\theta, r\sin\theta),
\]
where $G_m(x,y)$ is the homogeneous polynomial of degree $m$ defined by its Fourier expansion as
\[
 G_m(\cos\theta, \sin\theta):= -m\sum_{k=-m}^{k=m} \frac{c_{k}}{ki+ma}e^{k\theta i},
\]
and $G(\theta)=G_m(\cos\theta, \sin\theta)$ satisfies
$G'+maG+mR_m(\cos\theta, \sin\theta)=0$.

By using the above expression we get that the only algebraic
solution of system~(\ref{eq:8}) is the one that satisfies
\[
r_0^{-m}+m\sum\limits_{k=-m}^{k=m}\frac{c_{k}}{ki+ma}=0.
\]
Moreover, it is easy to check that the cofactor of this algebraic
solution, $H(x,y)=G_m(x,y)-1$, is $K(x,y)=mR_m(x,y)$, see also the proof of Lemma~\ref{le:2.7}

Notice that a necessary and sufficient condition for the
existence of a real algebraic solution is that
\[
\sum\limits_{k=-m}^{k=m}\frac{c_{k}}{ki+ma}<0.
\]

Finally, the limit cycle exists when $G_m(\cos\theta, \sin\theta)>0$ for $\theta\in[0, 2\pi].$ This only can happen when $m$ is even, see also \cite{GasTor2004}.
\end{proof}

\section{The examples}\label{se:4}
In this section we will prove that the system given in
Theorem~\ref{th:1.1} has an explicit limit cycle. That the
limit cycle is not algebraic will be proved in Section~\ref{se:6}.

The system~(\ref{eq:2}) is
\[
\left\{
\begin{array}{ccl}
\dot x&=&-(x-y)(x^2-xy+y^2)+x(2x^4+2x^2y^2+y^4), \\
\dot y&=&-(x+y)(2x^2-xy+2y^2)+y(2x^4+2x^2y^2+y^4),
\end{array}
\right.
\]
and in polar coordinates it can be written as
\[
\left\{
\begin{array}{ccl}
\dot r&=&(\cos^4 \theta +1)r^5+(\cos^2\theta -2)r^3, \\
\dot \theta&=&-(\cos^2\theta+1)r^2.
\end{array}
\right.
\]

By following the same steps than in the proof of
Theorem~\ref{th:2.1}, we introduce the change of variables
$r=1/\sqrt{\rho},$ obtaining:
\[
\rho'=2\frac{\cos^2\theta-2}{\cos^2\theta+1}\rho+2\frac{\cos^4\theta
+ 1 } { \cos^2\theta+1}.
\]
The solution satisfying that $\rho=\rho_0>0$ when $\theta=0$ is
\[
\rho(\theta;\rho_0) = e^{-3\psi(\theta)+2\theta}
\left(\rho_0+2
\int_0^\theta{\frac{\cos^4(s)+1}{\cos^2(s)+1}
e^{3\psi(s)- 2s}}ds\right)>0,
\]
where
$\psi(\theta)=2\int_0^\theta \frac{1+\tan^2s}{2+\tan^2s}ds.$

The initial condition of the limit cycle is given by the equation
$\rho(2\pi)=\rho(0)=\rho^*_0$. Hence
\[
\rho^*_0=\frac{2e^{4\pi-6\sqrt{2}\pi}}{1-e^{4\pi-6\sqrt{2}\pi}}
\int_0^{2\pi}\frac{\cos^4(s)+1}{\cos^2(s)+1}
e^{3\psi(s)-2s}ds>0.
\]
This value can be can be computed numerically, giving
$\rho_0^*\approx 1.1990$. The intersection of the limit cycle with
the OX$^+$ axis is the point having $r^*_0=1/\sqrt{\rho^*_0}\approx
0.9132$.

Since the Poincar\'e return map is
$\Pi(\rho_0)\!=\!\rho(2\pi;\rho_0)$ we have
$\Pi'(\rho_0)\!=\!e^{(4-6\sqrt{2})\pi}<1$ for all $\rho_0$ and
$\dot{\theta}<0$, we get that the limit cycle of system~(\ref{eq:2})
is hyperbolic and unstable.

\section{A method for studying the existence of algebraic solutions}\label{se:5}

Let $F(x,y),$ $K(x,y),$ $X(x,y)$ and $Y(x,y)$ be real analytic
functions such that
\begin{equation}\label{eq:9}
\frac{\partial F}{\partial x}X(x,y)+\frac{\partial F}{\partial y}
Y(x,y)= K(x,y)F(x,y).
\end{equation}
 Thus it is clear that the set $\{(x,y)\in
\mathbb{R}^2\, :\, F(x,y)=0\}$ is formed by solutions of the system
\begin{equation}\label{eq:10}
\left\{
\begin{array}{ccl}
\dot x&=&X(x,y), \\
\dot y&=&Y(x,y).
\end{array}
\right.
\end{equation}
Fixed an analytic solution of (\ref{eq:10}) of the form
$y=\alpha(x)$, we can consider the following Taylor expansions in
$z$,
\begin{eqnarray*}
F(x, z+\alpha(x))&=&F_0(x)+zF_1(x)+z^2F_2(x)+\ldots, \\
K(x, z+\alpha(x))&=&K_0(x)+zK_1(x)+z^2K_2(x)+\ldots, \\
X(x, z+\alpha(x))&=&X_0(x)+zX_1(x)+z^2X_2(x)+\ldots, \\
Y(x, z+\alpha(x))&=&Y_0(x)+zY_1(x)+z^2Y_2(x)+\ldots\\
\end{eqnarray*}
Notice that $\alpha'(x)=Y_0(x)/X_0(x)$. Then equation~(\ref{eq:9}) can be written as
\[
\sum\limits_{k=0}^{\infty}
\left(
\sum\limits_{i=0}^{k}
\Bigl(
X_{k-i}(x)F_i'(x)+\bigl(iY_{k-i+1}-i\alpha'(x)X_{k-i+1}-K_{k-i}\bigr)F_i(x)
\Bigr)
\right)
z^k=0.
\]
The functions $F_k(x)$ can be obtained
recurrently from the above relation by solving the linear differential equations in
$F_k(x)$, obtained vanishing each coefficient in $z^k.$

In particular, for $k=0$ and $k=1$, we get
\[
X_0(x)F_0'(x)-K_0(x)F_0(x)=0,
\]
\[
X_0(x)F_1'(x)+(Y_1(x)-\alpha'(x)X_1(x)-K_0(x))F_1(x)+X_1(x)F_0'(x)-K_1(x)F_0(x)=0.
\]
We obtain $F_0(x)=C_0\exp(\int_0^x \frac{K_0(s)}{X_0(s)}ds)$, where $C_0$ is an arbitrary constant, and similarly we could get $F_1(x)$.

When $\alpha(x)$ is a rational function and $F(x,y)$, $K(x,y)$,
$X(x,y)$ and $Y(x,y)$ are polynomials, with real or complex
coefficients, the linear differential equations for each $F_k(x)$
described in the above algorithm give us a collection of necessary
conditions for the existence of an algebraic solution $F(x,y)$.
The conditions are that, for each $k$, \textit{the functions
$F_k(x)$ must be polynomials}. For instance, for $k=0$, the first
necessary condition is that the primitive of the rational function
\[
\frac{K_0(x)}{X_0(x)}=\frac{K(x, \alpha(x))}{X(x, \alpha(x))}
\]
must be a linear combination of logarithms of polynomials. Furthermore
the coefficients of the logarithms have to be natural numbers.

The necessary conditions obtained for the existence of algebraic
solutions restrict the possible cofactors of $F.$ These
restrictions give the key for searching the possible algebraic
solutions of system~(\ref{eq:10}), see Remark~\ref{re:5.2}
and Section~\ref{se:6}. As we will see, in our case we only need
to apply the described method for $k=0$ but we remark that in
other situations, by using it for bigger $k,$ it can give more
information about the existence or non existence of algebraic solutions.

\begin{remark}\label{re:5.1}
Notice that the above method can only be applied
when the candidate $F(x,y)$ to be an algebraic solution of
system~(\ref{eq:10}) does not contain the factor $y-\alpha(x).$
\end{remark}

\begin{remark}\label{re:5.2}
Assume that system~(\ref{eq:10}) is fixed and it is polynomial. Notice that the equation~(\ref{eq:9}) that gives the possible set of algebraic
solutions of system~(\ref{eq:10}) is equivalent to a set
of quadratic equations where the unknowns are the coefficients of
$F$ and the coefficients of $K.$ In general, it is very hard to
solve this system of equations, even by using algebraic
manipulators. On the other hand, the method developed in this
section imposes restrictions on the cofactor $K$ for the existence
of $F.$ Ideally, if $K$ is totally known, the system to be
solved is be linear and so the problem of knowing the existence
or not of algebraic solutions of a given degree would be a much easier task. In any
case, any information on $K$ makes the problem simpler.
\end{remark}

\section{A non-algebraic limit cycle}\label{se:6}

This section is devoted to prove that the limit cycle of system~(\ref{eq:2}) is not algebraic. Indeed, by using the method introduced in Section~\ref{se:5} we will prove that the only
algebraic solutions of the system are the ones given in Lemma~\ref{le:2.4}. The algebraic solution given by this lemma is $yP_n-xQ_n=(2x^2+y^2)(x^2+y^2).$ Concretely, the curves $x^2+y^2=0$ and $2x^2+y^2=0$ have the cofactors $2(2x^4+2x^2y^2+y^4-x^2-2y^2)$ and
$2(2x^4+2x^2y^2+y^4-x^2+xy-2y^2),$ respectively. These two curves coincide with the four complex lines $y=\pm i x$ and $y=\pm \sqrt{2}i x.$

As we will see, the first step ($k=0$) of the method developed in
Section~\ref{se:5} applied to each one of the four complex lines, $y=\pm i x$ and $y=\pm \sqrt{2}i x,$ will give enough restrictions to prove
that the only algebraic solutions of system~(\ref{eq:2}) are
the ones described above.

Assume that the differential system has a real or complex algebraic solution $F$ and that it does not contain none of the given four lines as a factor. By using Lemma~\ref{le:2.6} it is not restrictive to assume that $F$ is real and that its cofactor is an even function, \textit{i.e.}
$K(-x, -y)=K(x,y).$

Since the degree of the vector field (\ref{eq:2}) is 5 we know that the degree of $K(x,y)$ is at most 4. By the above restrictions on $K(x,y)$ and by using also Lemma~\ref{le:2.5}
we can write it as the real polynomial
\[
K(x,y)=a_{00}+a_{20}x^2+a_{11}xy+a_{02}y^2+\ell(2x^4+2x^2y^2+y^4),
\]
where $\ell$ is the degree of the corresponding algebraic curve $F(x,y)=0$.

We apply the first step of our method, \textit{i.e.} we take $k=0.$ By considering the cases $\alpha(x)=\pm ix$ we obtain
\begin{eqnarray*}
 & & \int\frac{K(x, \alpha(x))}{X(x, \alpha(x))}dx=\int\frac{K(x, \pm ix)}{X(x, \pm ix)}dx=\int\frac{K_0(x)}{X_0(x)}dx=\\
  \end{eqnarray*}
  \begin{eqnarray*}
 & &\frac{a_{00}(-1\pm i)}{4x^2}+\frac{1}{2}(a_{20}+a_{11}-a_{02}\pm i(-a_{20}+ a_{11}+a_{02}+ a_{00})) \log(x)+\\
 & &\frac{1}{8}(-a_{20}-a_{11}+a_{02}+2\ell\pm i( a_{20}- a_{11}- a_{02}- a_{00})) \log(2+2 x^2+x^4)+\\
 & &\frac{1}{4}(a_{20}-a_{11}-a_{02}-a_{00}\pm i(a_{20}+a_{11}-a_{02}-2\ell)) \arctan(x^2+1).\\
 \end{eqnarray*}
By forcing
$F(x, \alpha(x))\!=\!F(x, \pm ix)\!=\!F_0(x)\!=\!C_0\exp\left(\int K_0(x, \pm ix)/X_0(x, \pm ix)\, dx\right)$
to be a polynomial, with $C_0$ an arbitrary constant, we obtain a first set of necessary conditions:
\[
\left\{\begin{array}{lcr}
 a_{20}-a_{11}-a_{02}-a_{00}&=&0, \\
 a_{20}+a_{11}-a_{02}-2\ell &=&0, \\
 a_{00}&=&0.\\
 \end{array}
 \right.
 \]
The same computations can be done for the other pair of
algebraic solutions, $y=\pm\sqrt{2}ix,$ that is
 \begin{eqnarray*}
 & & \int\frac{K(x, \alpha(x))}{X(x, \alpha(x))}dx=\int\frac{K(x, \pm i\sqrt{2}x)}
 {X(x, \pm i\sqrt{2} x)}dx=\int\frac{K_0(x)}{X_0(x)}dx=\\
 & & -\frac{a_{00}}{6x^2}+\frac{1}{9}(3 a_{20}-6a_{02}-2a_{00}\pm 3 \sqrt{2} i a_{11} )\log(x)+ \\
 & & \frac{1}{18}(-3 a_{20}+6 a_{02}+9\ell+2 a_{00}\mp 3\sqrt{2} i a_{11} ) \log(3+2 x^2).\\
 \end{eqnarray*}
As in the previous case, we obtain a second set of necessary conditions:
\[
\left\{\begin{array}{lcr}
 a_{11}&=&0, \\
 a_{00}&=&0.\\
 \end{array}
 \right.
 \]
Collecting all the obtained equations we get that the degree of
the invariant algebraic curve is $\ell=0$, or in other words that
such a curve does not exist.

\vspace{5mm}

\noindent\textbf{Proof of Theorem~\ref{th:1.1}}. The proof
of the theorem follows from the results of Sections~\ref{se:2},
\ref{se:4} and \ref{se:6}.

\begin{remark}\label{re:6.1}
>From \cite[Thm 1]{Sin1992}, if a planar system has an explicit
non-algebraic solution which is in the zero level set of a
Liouvillian function then it has a Darboux integrating factor and
therefore the whole system is integrable by quadratures. Notice
that this is the case for system (\ref{eq:2}): it has a
non-algebraic Liouvillian limit cycle and it can be transformed
into a Bernoulli equation. Consequently, if we would like to have
an explicit non-algebraic limit cycle for a planar system, which
is not integrable by quadratures, we should look for a limit cycle
given by a non-Liouvillian function.
\end{remark}

\vspace{1cm}

\end{document}